\newcommand{\beq}{\begin{quote}}
\newcommand{\enq}{\end{quote}}
\newcommand{\be}{\begin{equation}}
\newcommand{\en}{\end{equation}}
\newcommand{\del}{\delta}

\newcommand{\eps}{\epsilon}

\documentclass[11pt]{amsart}
\usepackage{geometry}                
\usepackage{graphicx}
\usepackage{amssymb}
\usepackage{epstopdf}
\title{Barrow and Leibniz on  the fundamental theorem of the Calculus}
\date{Nov. 2010}
\begin{document}
\maketitle
\begin{abstract}
In 1693,  Gottfried Whilhelm Leibniz published in the {\it Acta Eruditorum}  a {\it geometrical}  proof of the fundamental theorem of the calculus. 
During  his notorious  dispute with Isaac Newton on the
development of the calculus,  Leibniz  denied  any indebtedness to the work of  Isaac  Barrow. But  it is shown here, that 
his   geometrical proof of this theorem closely resembles  
Barrow's proof  in  Proposition 11, Lecture 10,  of his  {\it Lectiones Geometricae},
published in 1670.

 \end{abstract}

\section{Introduction}

At the height of his priority  dispute with Newton concerning the invention of the calculus,
Leibniz wrote an account,  {\it Historia et Origo Calculi Differentialis},   describing the 
 contributions  by  seventeenth  century mathematicians that led him to his own development of the calculus (Child 1920, 22).  In this account,
 Isaac Barrow  is not mentioned at all,  and in several occasions Leibniz  denied any indebtedness to his work,  particularly during his
 notorious priority dispute with Isaac Newton\footnote{
 Earlier, however, Leibniz did refer  Barrow's work.  In his 1686 article  in the {\it Acta Eruditorum},  ``On a  deeply hidden
 geometry and the analysis of indivisibles and infinities", Leibniz referred  to Barrow in connection with
 a geometrical theorem that appeared  in Barrow's 
 Geometrical Lectures,  and proceeded  to give a proof of this theorem by his analytic method (Struik 1696, 281). Also, on
 Nov. 1, 1675, Leibniz wrote: ``Most of the theorems of the geometry of indivisibles which are to be found in
 the works of Cavaliere, Vincent, Wallis, Gregory and Barrow, are immediately evident from the calculus" (Child 1920, 87).
 I am indebted to an anonymous reviewer for calling my attention to these two references.
 }.
 But in  Barrow's {\it Lectiones Geometricae} (hereafter cited as {\it Geometrical Lectures}),  which Leibniz had obtained  during a visit to London in 1673,  the 
 concepts of the differential and integral calculus are discussed in {\it geometrical} form,  and a  rigorous mathematical proof 
 is given  of the fundamental theorem of the calculus\footnote
 {This  theorem establishes  the inverse relation between integration and differentiation.  In geometrical terms, it  equates   the  subtangent 
 of a curve $\alpha$  that gives  the area enclosed by   a given  curve $\beta$, to the ratio of the ordinates
  of these two curves (see section I).
Introducing Cartesian coordinates  $x,z$ for  $\alpha $, and $x,y$ for $\beta$, where $x$ is the common abscissa, 
the subtangent $t$  of $\alpha$ is $t=z dx/dz$, and expressed  in analytic form,
 Barrow's theorem establishes  the relation
 $ t=z/y$. Hence, Barrow's theorem is equivalent to   the relation  $dz/dx =y$ for 
 the fundamental theorem of the calculus.} \footnote
 {It should also be pointed out that a geometrical proof of the fundamental theorem of the calculus similar to Barrow's 
 was given by the Scottish mathematician James Gregory, which he  published   as 
 Prop. VI in his book {\it Geometriae pars universalis} (Padua 1668) (Baron  1969, 232).
 }.
  According to J. M. Child,  ``a Calculus may be of two kinds:
 \beq
 
 i) An {\it analytic} calculus, properly so called, that is, a set of algebraical working rules (with their proofs), with
 which differentiations of known functions of a {\it  dependent} variable, of products, of quotients, etc., can be 
 carried out; together with the full recognition that differentiation and integration are inverse operations, to enable
 integration from first principles to be avoided . . .
 
 ii) A {\it geometrical} calculus equivalent embodying the same principles and methods; this would be the more
 perfect if the construction for tangents and areas could be immediately  translated into algebraic form, if it
 where so desired.
 \enq
 Child concludes
 \beq
 Between these two there is,  in my opinion,  not a pin to choose theoretically; it  is a mere matter of practical utility that
 set the first type in front of the second; whereas the balance of rigour, without modern considerations, is
 all on the side of the second ( Child 1930, 296).
 \enq
More recently,  a distinguished mathematician, Otto Toeplitz,  wrote that `` Barrow was in possesion of most of the rules of 
 differentiation, that he could treat many inverse tangent problem (indefinite integrals), and that in 1667 he discovered
 and gave an admirable proof of  the fundamental theorem - that is  the relation [ of the inverse tangent] to the definite
 integral (Toeplitz 1963,  128).
  
Leibniz's  original work concerned the {\it analytic} calculus,  and he claimed to have  read the relevant  sections of Barrow's lectures on 
the {\it geometrical} calculus  only several years later, after he 
 had independently made his own discoveries.  In a letter to Johann Bernoulli  written in 1703 , he attributed his initial inspiration
 to a ``characteristic triangle"   he had found in Pascal's {\it Trait\'{e} des sinus du quart the cercle.}\footnote
 {  Leibniz said that on the
 reading  of this example in Pascal a light suddenly burst upon him, and  that  he then realized
 what Pascal had not  - that the determination  of a tangent to a curve depend on the ratio
 of the differences in the ordinates and abscissas, as these became infinitesimally small,
 and that the quadrature depended upon the sum  of ordinates or inifinitely thin rectangles for 
 infinitesimal intervals on the axis. Moreover, the operations of summing and of finding 
 differences where mutually inverse (Boyer 1949, 203).}
 After examining Leibniz's original manuscript, and his
 copy of  Barrow's  {\it Geometrical Lectures},  which
 were not available to Child,   D. Mahnke (Manhke 1926) concluded that Leibniz had read only the beginning of 
 this book,  and that his  calculus discoveries were made independently 
 of Barrow's work.
 Later,    J. E. Hofmann (Hofmann 1974),  likewise concluded  that this independence
 is confirmed  from Leibniz's  early  manuscripts and correspondence at the time (Hofmann 1974).  

 In his early mathematical studies,  Leibniz  considered number sequences  and realized that the operations
 associated with certain sums and  differences of these sequences had a reciprocal relation. Then, by  approximating a curve
 by polygons,  these  sums correspond to the area bounded by the  curve, while the differences correspond 
 to its  tangent, (Bos 1973; Bos 1986, 103),  which Leibniz indicated \footnote
 { See {\it Historia et origo calculi differentialis} (Child 1920, 31-34)}, 
 lead to his insight of the reciprocal relation between areas and tangents .  For example, 
 Leibniz's  June 11, 1677 letter addressed to Oldenburg for Newton, in reply to  Newton's  October 24, 1676 letter to Leibniz ({\it  Epistola Posterior}), clearly shows that by this time he understood the fundamental theorem of the calculus, and its usefulness
  to evaluate algebraic integrals (Newton 1960, 221; Guicciardini  209, 360).
 But  Child has given considerable 
 circumstantial evidence, based on the reproductions of some of Leibniz's  manuscripts,  that Barrow's {\it Geometrical Lectures}  also influenced some of Leibniz's work\footnote
 {In  the  preface to his 1916   translation of  Barrow's {\it Geometrical Lectures},  Child concluded that
\beq
``Isaac Barrow was the first inventor of the Infinitesimal Calculus; Newton got the main idea of it from 
Barrow by personal communication;  and Leibniz also was in some measure indebted to  Barrow's work,
obtaining confirmation of his own original ideas, and suggestions for their further development, from the
copy of Barrow's  book that he purchased in 1673."  (Child 1916,  7)
\enq
}
(Child 1920), and he did not concur with Mahnke's conclusion\footnote
{  In 1930, in an article written partly in response to Mahnke's article, Child concluded:
\beq
``But if he [Leibniz]  had never seen Barrow,  I very much doubt  if he (or even Newton) would
have invented the analytic calculus, and completed it in their lifetimes." (Child 1930, 307)
\enq
}.  
To give one example, in  the first publication of his integral calculus (Leibniz 1686),  Leibniz gave an analytic derivation
of Barrow's geometrical  proof in Prop. 1, Lecture 11, (Child 1916, 125),  that the area bounded by a curve with ordinates equal to the subnormals of a {\it given} curve
is equal to half  of the square of the final ordinate of  the original curve\footnote{
Let $y(x)$ be the ordinate of a given curve, where $y(0)=0$.  The subnormal is $n(x)=y dy/dx$, and in accordance
with Leibniz's rules for integration, $\int^{x}_0 n(x') dx'=\int^{y}_0 y' (dy' /dx') dx'=\int^{y}_0 y'dy' =(1/2)y^2$.
This example illustrates the greater  simplicity of Leibniz's analytic derivation, based on his  suggestive notation,
which treats differentials like $dx'$ as algebraic quantities,  compared to Barrow's 
geometrical formulation, see Appendix A}.  In this publication, 
Leibniz admitted that he was familiar 
with Barrow's theorem\footnote {Child calls attention to the fact that  the ``characteristic triangle", which Leibniz claimed to have learned from
Pascal, appears also in Barrow's diagram in a form very similar to that described  by Leibniz (Child 1920, 16), and concluded, rather dramatically, 
that ``such evidence  as that would be enough  to hang a man, even in an English criminal court".
},
and notes found in the margin of his copy of Barrow's  {\it Lecciones Geometricae} indicate that Leibniz  had read it  sometimes between 1674 and 1676 (Leibniz 2008, 301).  Barrow's  proposition appears in one of he last lectures in his book, which suggests that by this time he had read also most of the other theorems in  this book,  particularly  Prop. 11, Lecture 10 which contains  Barrow's geometrical proof of the fundamental
theorem of the calculus.

 But in a 1694 letter to the Marquis de l'Hospital,  Leibniz  wrote:
 \beq
 I recognize that M. Barrow has advanced considerably, but I can assure you, Sir, that I have derived no assistance  from him for my methods
 (pour mes methodes) (Child 1920, 220),
 \enq
 and insisted that
 \beq
 by the use of the ``characteristic triangle" . . . I thus found  as it were in the twinkling of an eyelid nearly all the theorems that I afterward found
 in the works of Barrow and Gregory (Child 1920, 221).
 \enq
 Later, in an intended  postscript\footnote
 { This postcript was in a draft,  but not included in the letter sent to Jacob Bernoulli.   Child  erroneously stated  that this letter
 was sent to his brother, Johann Bernoulli. I am indebted to a reviewer for this  correction.
 }
 to a letter from Berlin to Jacob Bernoulli, dated April 1703,  he wrote,
 \beq
 Perhaps you will think it small-minded of me that I should be irritated with you, your brother [Johann Bernoulli], or anyone else,
 it you should have perceived the opportunities for obligation to Barrow, which it was not necessary for me, his contemporary
 in these discoveries\footnote
 { In this letter  Leibniz contradicted himself,  because earlier  he had  claimed that his discoveries had occurred several years after Barrow's book
 had appeared.}
  to have obtained from him. (Child 1920,  11)
 \enq
Apparently this comment was motivated  by an article in  the {\it Acta Eruditorum} of
January 1691,  where Jacob Bernoulli  had emphasized the similarities of the methods of  Barrow  with those of Leibniz, and
argued that anyone who was familliar with the former could hardly avoid recognizing the latter,
\beq
Yet to speak frankly, whoever has understood Barrrow's . . . will hardly fail to know the other discoveries of Mr. Leibniz considering that they were based on that earlier discovery, and do not differ from them, except perhaps in the notation of the differentials and in some abridgment of the operation of it.(Feingold 1993,  325)
\enq 

More recently, M. S.  Mahoney,  presumably referring  to Child's work,  remarked   that
``Barrow seems to have acquired significant historical importance only at the turn of the twentieth century, when historians
revived his reputation on two grounds: as a forerunner of the calculus and as a source of Newton's mathematics." 
But, he continued,  ``beginning
in the 1960's two lines of historical inquiry began to cast doubt on this consensus" , and  he concluded that
Barrow was `` competent and well  informed, but not particularly original" (Mahoney 1990,  180, 240).  This sentiment echoes Whiteside's  assessment
that  ``Barrow remains  . . . only a  thoroughly competent university don whose real importance lies more in his coordinating
of available  knowledge for future use rather than in introducing  new concepts"   (Whiteside 1961, 289). In contrast to these derogatory comments,
Toeplitz concluded, in conformance with Child's evaluation, that ``in a very large measure Barrow  is indeed the real discoverer  [of the fundamental theorem of the calculus] - insofar as an individual can ever be given credit within a course of development such as we have tried to trace here" (Toeplitz  1963, 98).
 Whiteside claimed  
that Barrow's  proof of the fundamental  theorem of the calculus is only a ``neat amendment of [James] Gregory's generalization of Neil's rectification method,"  where  Barrow ``merely replaced an  element of arc length by the ordinate of a curve.''  
But M. Feingold has pointed  out  that Whiteside's claim that Barrow borrowed  results from Gregory  is not consistent with the fact that Barrow's manuscript  was virtually finished and already in the hands of John
Collins, who shepherded it  trough  the publication process  by  the time that Barrow received Gregory's book\footnote{
In Prop. 10, Lecture 11, however,  Barrow remarks: ``This extremely useful theorem is due to that most learned man, Gregory from Abardeen"
indicating that he was familiar at least with one  of Gregory's geometrical proofs.  Interestingly, as Child pointed out (Child 1920, 140),  Leibniz also refers to this 
theorem as,  `` an elegant theorem due to Gregory". }(Feingold 1993, 333).

This article  is confined primarily to Leibniz  {\it geometrical}  discussion of the fundamental theorem of the
 calculus as it  appeared in  the 1693 {\it Acta Eruditorum},  and  to a detailed  comparison
 with Barrow's  two  proofs of this theorem in  Proposition 11, Lecture 10, and in Prop. 19, Lecture 11, published
in his {\it Geometrical Lectures} 23 years earlier.  In particular,  Leibniz's early development of the analytic form of the calculus, and his
introduction of the  suggestive  notation that is in use up to the present time,  will not be discussed here\footnote{
 Leibniz's  development of the analytic form of the calculus in the period 1673-1676 has been discussed in great detail in the past.
For an excellent account, see  Christoph J. Scriba's article  {\it The Inverse method of Tangents:  A Dialogue between Leibniz and Newton (1675-1677)}, in
Archive for History of Exact Sciences $\bf 2$, 1963, 112, 137.
}, and   likewise, 
the  related work  of Isaac Newton,  which  recently has been covered in great detail   by N. Guicciardini  in his book {\it Isaac Newton on Mathematical Certainty and Method},
that  contains also  an excellent set of references to earlier work on this subject  (Guicciardini 2009)
In  the next section,   Barrow's geometrical proofs of these two propositions is discussed, and 
in Section 3  the  corresponding demonstration by
Leibniz is given, showing its  close  similarity to  Barrow's  work.   This resemblance 
has not been discussed  before,  although it was  noticed recently by L. Giacardi (Giacardi  1995), and it  supports  Child's thesis that  some of Leibniz's contributions to the development of  the calculus were 
influenced by Barrow's work. Section 4 gives a description of an ingenious mechanical device  invented by Leibniz to
obtain the area bounded by a given curve, and 
Section 5  contains a  brief summary and  conclusions.  In Appendix A  a detailed  discussion is given of Barrow's proposition
1, Lecture 11, that illustrates his application of the characteristic triangle,  which Leibniz  obtained sometimes between 1674 and 1676. 
Finally, in Appendix B  it is  pointed out  that  significant errors 
occur in the  reproduction of Leibniz's diagram associated with the fundamental theorem as it appears in Gerhard's edition (Leibniz 1693),  in Struik's English translation
(Struik 1969, 282),  and more recently in L. Giacardi's book (Giacardi 1995),  that apparently have remained  unnoticed in the past.


\section{  Barrow's geometrical proof of the fundamental theorem of the calculus}

 \begin{figure}[htbp] 
   \centering
   \includegraphics[width=4in] {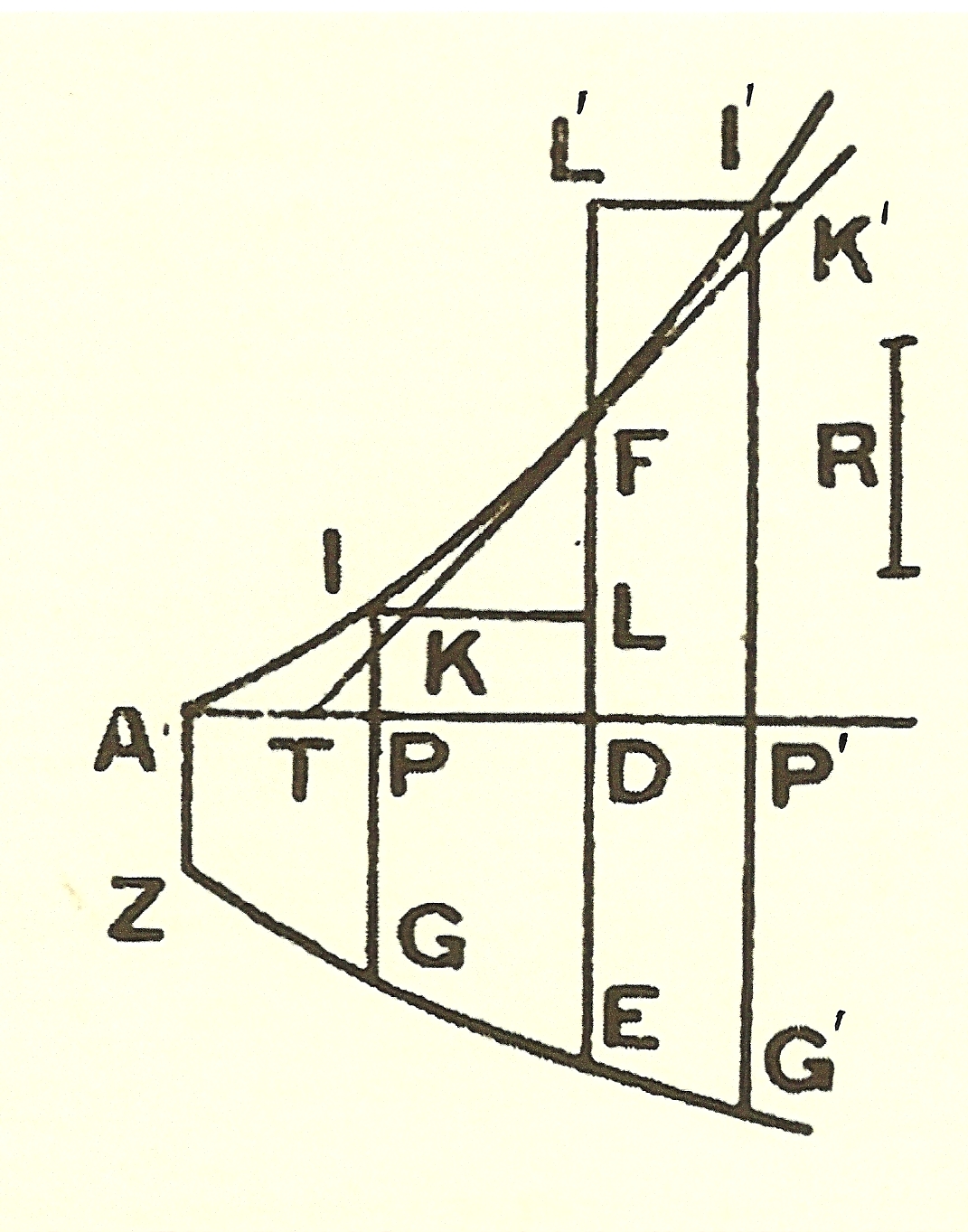} 
   \caption{Barrow's  diagram for Proposition 11, Lecture 10.  For clarity,
    prime superscripts have been added to labels that appear repeated  in the original diagram. }
   \label{fig:example}
\end{figure}

In lectures 10 and 11 of  his {\it Geometrical Lectures},   Barrow gave two related geometrical proofs of the fundamental theorem
of the calculus, (Child 1916, 116,135; Struik 1969, 256; Mahoney 1990, 223,232; Folkerts 2001; Guicciardini 2009, 175).  Underscoring
the value  of this theorem,  he started his discussion with the remark
\beq
I add one or two theorems which it will be seen are of great generality, and not lightly to be passed over\footnote{
Nevertheless, Mahoney asserted that ``what in substance becomes part of the fundamental theorem of the calculus
is clearly not fundamental to Barrow" (Mahoney 1990, 236).
}
\enq
Referring to his diagram  in Prop. 11, Lecture 10, shown\footnote
{ Prime superscripts have been added  to  the  labels $I,K,L,P$ and $G$ that  appeared repeated 
on the right side of the original diagram}
  in  Fig. 1, the area bounded by a given curve $ZGEG'$  is 
described by  another   curve $AIFI'$ with a common abscissa $APDP'$.   Barrow  then proved that the    subtangent $DT$ of $AIFI' $ at $F$
is  proportional to the ratio of the  ordinates of  the curves $AIFI'$ and $ZGEG'$, at a corresponding value $AD$ of the abscissa, i.e. $DT =R(DF/DE)$, where $R$ is an arbitrary parameter with dimensions of length\footnote
{ In Cartesian notation $AD=x, DE=y, DF=z$,  Barrow's theorem states that the subtangent $DT=t$ of
$AIFI' $ at $x$ satisfies the property $t=R(z/y)$. Treating
$dx=IL$ and $dz=FL$ as differentials,  i.e.,  setting  $KL\approx IL$,  by the similarity of triangles $TDF$ and $KLF$ we have 
 $t/z\approx dx/dz$. Hence, setting $R=1$, Barrows theorem corresponds to the familiar
analytic  form  of the fundamental theorem of the calculus  that, in the limit that $dx \rightarrow 0$,  $dz/dx=y$.}.
Even by modern standards,  Barrow's  proof  is  mathematically rigorous,  and  it does not depend on any  assumptions about differentials or  infinitesimals which were not well understood at the time.\footnote
{  Child wrote   that `` . . . it only remains to remark on the fact that the theorem of Prop. 11 is a rigorous proof  that differentiation
and integration are inverse operations,  where integration is defined as a summation." (Child 1916, 124) In Prop. 11, Lecture 10,
  Barrow  did not  express the area or integral  as 
 a summation, but  in Prop. 19, Lecture 11, he expressed integration as a  summation. 
} 
 In Proposition 19, Lecture 11,  Barrow gave a related  geometrical proof of his theorem  that was based  
 on an expression for the subtangent and  for the area of curves in terms of differential quantities that were  well known  
to mathematicians in the seventeenth century (Boyer, 1959).

Barrow started Prop. 11  with a description of his diagram,  shown in Fig. 1, as follows: 
 \beq
Let $ZGE$ be any curve of which the axis is $AD$, and let ordinates applied to this axis $AZ, PG, DE$ continually increase 
from the initial coordinate $AZ$; and also let $AIF$ be a line [another curve] such that, if any line $EDF$ is drawn perpendicular
to $AD$ cutting the curves  in the points $E, F$, and $AD$ in $D$ the rectangle contained  by $DF$ and a given length $R$
is equal to the intercepted space $ADEZ$;'' (Child 1916, 117)
\enq
In other word,  the ordinate $DF$  that determines   the curve $AIFI'$  is  equal to the area\footnote 
{ Barrow's word  for  area is ``space."  } $ADEZ$ bounded by the curve $ZGE$, the abcissa
$AD$ and the ordinates $AZ$ and $DE$.  Because the quantity associated with an area has the dimensions of length squared,  Barrow introduced  a parameter $R$ with units of length,  and set
\be
 DF= area ( ADEZ)/ R.
\en
  Barrow concluded, 
\beq
also let $DE/DF=R/DT$, and join $DT$. Then $TF$ will touch the curve $AIF$(Child, 1916, 117),
\enq
which meant  that   $TF$  is a line tangent to the curve $AIF$ at $F$.
Since it is  $DE$  and $DF$ that are  determined at a given valued of the abscissa $AD$,  $DT$ is {\it defined} by the relation
\be
\label{subt1}
DT=R\frac{DF}{DE},
\en
where $R\cdot DF=$ area $ADEZ$.  Hence,  Barrow's  proof  
consisted  in showing that   $DT$ is the subtangent  of the curve $AIFI'$ at $F$.  In   Section 3 , we   shall see that  Leibniz named this relation  the ``tangency relation"  without, however, crediting it to Barrow's original work. 

Before following Barrow's proof further,  it is   interesting  to speculate  how he  may have discovered that his relation
 for  $DT$, Eq. \ref{subt1}, corresponds to the {\it subtangent} at $F$ of the curve  $AIFI'$ in the limit that $DP$ becomes vanishingly small. 
At the end of  lecture $10$,  Barrow gave  an argument, originally  due to  Fermat, that 
when $I$ and  $F$ are close to each other, the triangles $TDF$ and $ILF$ are approximately similar, and therefore,\footnote
{ The symbol $\approx$ is introduced  to  indicate  that this relation is only approximately valid for finite differentials $IL,FL$.}

\be
\label{subt}
\frac{DT}{DF} \approx \frac{IL}{FL}.
\en
Leibniz named a triangle similar to  $ILF$ the ``characteristic triangle",  but   attributed it to
Pascal.

Since by  definition  of the curve $AIFI'$, 
\be
FL=DF-IP=\frac{area (PDEG)}{R},
\en
and
\be
area(PDEG)\approx DE\cdot PD
\en
where $PD=IL$,  we have 
\be
FL\approx \frac{DE \cdot IL}{R}.
\en
Substituting  this approximation for $FL$ in Eq. \ref{subt} leads to the  relation
\be
\frac{DT}{DF}\approx \frac{R}{DE}.
\en
In the limit that $I$ approaches $F$,  and $IL$ becomes vanishingly small, the approximate sign $\approx$ in this equation becomes
 an  equality, and  this expression  becomes  Barrow's   relation for the subtangent $DT$, Eq.  \ref{subt1}.

In the next step, Barrow gave   a rigorous proof for this relation, without appealing to differentials.
For  the
case that $DE$  increases with increasing value of $AD$, Barrow applied the inequality $DE \cdot PD > area (PDEG) $,  to show that for any finite value of $IL$, $LK< LI$   when $I$ is  on the left hand side of $F$,  and when $I$ is on the right hand side $F$ then  $LK>LI$.  Since $K$ lies on
the tangent line, and the curve $AIF$ is convex,  this result implies that the line $TF$,  that by definition ``touches"  the curve at $F$  
does not cross it at any other point. Therefore $TF$  is the tangent line at $F$ of the curve $AIF$.  
\beq
For, if any point $I$ is taken in the line $AIF$ (first on the side of F towards A), and if through it $IG$ is drawn parallel to $AZ$ and $KL$
 is parallel to $AD$, cutting the given  line as shown in the figure, then
 
 $LF:LK=DF:DT=DE:R$
 or
  $R\cdot LF= LK\cdot DE$
 
 But, from the stated nature of the lines $DF, PK$, we have $R\cdot LF= area (PDEG)$;
 therefore $LK \cdot DE =area (PDEG) < PD \cdot DE$;  hence $LK < DP=LI$.

 Again, if the point  $I$ is taken on the other side of  $F$, and the same construction is made as before, plainly
 it can be easily shown that $LK>DP=LI$.
 
 From which it is quite clear that the whole of the line $TKFK$ lies with or below  the curve $AIF$. (Child 1916,117)
 \enq
 At the end of this presentation,  Barrow  indicated  that  when $DE$  decreases with increasing values of $AD$ 
 \beq
  the same conclusion
 is attained by similar arguments; only one distinction occurs,  namely, in this case, contrary to the
 other, the curve $AIF$ is concave to the axis AD.(Child 1916, 118)
 \enq
 He concluded his proposition with the corollary 
 \be
 DE \cdot DT= R\cdot DF = area (ADEZ),
 \en  
which follows from his relation, Eq. \ref{subt1}.

At the end of lecture 10,  Barrow  considered the application of his geometrical
 result to the case that the curve $AIFI'$ is determined by an algebraic relation between the abscissa
$AD$ and the ordinate $DF$. Following Fermat, he set  $a=IL$ and $e=FL$, and  substituting $AD$  for the abscissa,  and $DF$ for the
ordinate,  he showed  by a set of three rules how to obtain  an expression for the ratio $e/a$ in terms of $AD$ and $DF$.
These rules amount to keeping only terms that are linear 
in $e$ and $a$ in a power series expansion of this algebraic relation\footnote
{ Barrow's  rules correspond to the modern definition of the derivative. Let  $z'=z-e$ be the ordinate $DF=z $ at a  value
$x'=AD-a$ of the abscissa $AD=x$.  Then $e/a=(z'-z)/a$,  and Barrow's rule is to expand $z'-z$ in powers of 
$a$ and neglect any term on the right hand side of this relation that depends on $a$.
 }.        

It should be pointed out  that  in his proof, Barrow did not have to specify how to evaluate
the area $(ADEZ)$, e.g.  by the sum of differential rectangles  indicated by  Leibniz's notation.
  Hence,  in  Prop. 11 Barrow gave  a rigorous proof
of the fundamental theorem of the calculus which in 
Cartesian and differential coordinates can be expressed in the following form: 
\beq
Given a curve $x,y$, there exists another  curve
$x,z$,  where  $z$ is the area of the region bounded by the given curve and its  coordinate lines, that 
has  the property
 that its derivative $dz/dx$  is proportional to $y$.
\enq

Along similar lines, in Proposition 19, Lecture 11, Barrow gave  another   proof of the fundamental theorem of the
calculus where he explicitly made use of differentials. For this proposition, Barrow introduced a somewhat different diagram,   but 
 since this is not necessary,  I will present his discussion  by using  his diagram, Fig. 1, from Prop. 11, Lecture 10. This  procedure will also facilitate comparisons
in the next section between Leibniz's and Barrow's proofs,  because Leibniz used essentially Barrow's  diagram 
in Prop. 11, Lecture 10. 
In Prop. 19, Lecture 11,  the order in which   Barrow  constructed  the two curves  also was  altered : he assumed that  $AFI$
is the given curve, and obtained  the  associated curve $ZEG$  by implementing 
relation,  Eq. \ref {subt1},   for
the ordinate  $DE$ of this curve in terms  of  the
subtangent $DT$ and ordinate $DF$ of the given curve$AFI$,  which previously had been derived quantities.  In the next section it will be shown that Leibniz's
geometrical proof  of the fundamental theorem of the calculus  is essentially the same as the proof that
 Barrow gave  in  Prop, 19, Lecture 11.

Barrow started Prop.19, Lecture   with the following  description of his geometrical construction, 
\beq
Again, let $AFI'$ be a curve  with  axis  $AP$ and let $PG$ be perpendicular to $AP$;  also let $ZEG'$ be another line [curve]
such that, when any point $F$ is taken on the curve $AFI'$ and through it are drawn $FT$ a tangent to the curve $AFI'$ and $FDE$
parallel to $PG$, cutting $ZEG'$ in $E$ and $AP$ in $D$, and $R$ is a line of  given length, $DT:FD=R:DE$. (Child 1916, 135) 
\enq
Here the ordinate $DE$ for the curve $ZEG'$ is determined by the Barrow's tangency relation, Eq. \ref{subt1},  
 previously established in Prop. 11, Lecture 10.  Barrow  formulated the fundamental theorem  as follows:

\beq
Then the space $AP'G'Z$
is equal to the rectangle contained by $R$ and $I'P$. (Child 1916, 135)
\enq

Here  Barrow's  proof made explicit application of differentials:

\beq
 For if $IF$ is taken to be an indefinitely small arc of the curve $AFI'$  [and $IL$ is drawn parallel to $AP$ cutting $FD$ at $L$]. . . then we have

$IL:FL:=TD:FD=R:DE$;
therefore $ IL\cdot DE= FL\cdot R$, and $PD\cdot DE=FL\cdot R$. (Child 1916, 135)
\enq
Setting $DE=y$, $PD=dx$ and $FL=dz$, Barrow's  relation takes the form
\be
dx\cdot  y=R\cdot dz.
\en

 Barrow applied the approximate similarity of the differential triangle $ILF$ and   the triangle $TDF$ associated
with the tangent line of $AFI'$ at $F$,  to equate 
the area of the differential rectangle $PD\cdot DE$  to the differential change $FL=FD-IP$ in the ordinate of $AFI'$.  Barrow concluded
\beq
Hence, since the sum of such rectangles as  $PD\cdot DE$ differs only in the least degree from the space $AP'G'Z$, and
the [sum of the]  rectangles $FL\cdot R$ from the rectangle $IP\cdot R$, the theorem is quite obvious (Child 1916, 135)
\enq
In Leibniz's notation this  ``sum of rectangles" is expressed in the form
\be
\int dx\cdot  y=R\cdot \int dz= R\cdot z,
\en
corresponding, in modern notation,  to the integral  relation for the area
bounded by a curve.  

Although Prop. 19  has been described as the converse
of Prop 11,  this characterization misses the relevance of this proposition  to complete the
formulation of the fundamental theorem of the calculus. 
For finite differentials, the product $PD\cdot DE$ in Barrow's proof is  larger than the area of the region $PDEG$; 
therefore the sum of the areas
of these rectangles  only gives an upper bound to the area of the region $AP'G'Z$.  
In the second  appendix to Lecture 12, Barrow  gave a geometrical proof that the area bounded by a concave curve,  is obtained by the 
`` indefinite" sum of either  circumscribed or inscribed rectangles. By indicating  that these two sums would be  the same, Barrow {\it outlined}
a geometrical proof for the existence of the integral of such
a curve, but his proof  was not completed until 17 years later by   Newton,  in Section 1, Lemma 2 of his {\it Principia}.
Finally, in 1854  Riemann extended the  method of Barrow and Newton to describe the necessary properties  of a function for which the
concept of an integral as an indefinite sum of rectangles is justified.
\footnote
{ In his 1854 Habilitationsschrift,  published posthumously by Dedekind, Riemann  introduced a bound similar to Barrow's  to establish the conditions for the existence of the
integral of a function $f(x)$. He wrote:

``First: What is one to understand by
$\int_a^b f(x) dx $? In order to fix this relation, we take between $a$ and $b$ a series of values $x_1,x_2, . . . x_{n-1}$, and
describe the short intervals $x_1-a$ by $\del_1$, $x_2-x_1$ by $\del_2$, . . ., $b-x_{n-1}$ by $ \del_n$.
Hence, the value of the sum
$S=\del_1 f(a+\eps_1\del_1)+\del_2 f(x_1+\eps_2\del_2)+. . .+\del_nf(x_{n-1}+\eps_n\del_n)$
will depend on $\del$ and the magnitude of $\eps$ [ remark: $0\leq \eps \leq 1$]. Given the property that when the $\delta's$
become vanishingly small, the sum approaches a limit  $A$, then this limit corresponds to $\int_a^b f(x) dx $.
If it does not have this property, then $ \int_a^b f(x) dx$  does not have any meaning.

Under what conditions will a function $f(x)$  permit an integration,  and when will it not?
Next, we  consider the concept of an integral in a narrow sense,  that is, we examine the convergence
of the  sum $S$,  when the various values of $\delta$ become vanishingly small.  Indicating  the largest
oscillation of the function $f(x)$ between $a$ and $x_1$, that is the difference between its largest and
smallest value in this interval by $D_1$,  between $x_1$ and $x_2$ by $D_2$ ..., between
$x_{n-1}$ and $b$ by $D_n$; hence the sum
$\del_1 D_1+\del_2 D_2+. . . +\del_n D_n,$
with the largest value of $\delta$,  must become vanishingly small. Furthermore, we assume that as  long as the $\delta's $ remain
smaller than $d$,  then the   greatest value of this sum is $\Delta$;  $\Delta$ will then be a function of $d$ which decreases in
magnitude with $d$ and with its size unendlessly decreasing." (Riemann 1854)

Barrow and Newton divided the interval $[a,b]$  into ``equal parts'', i.e.    $\del_i=d$, and for the special case that the function
 $f(x)$ decreases monotonically with increasing $x$, the Riemann's sum is equal to $f(a)\cdot d$.

}


\section{Leibniz's geometrical proof of the fundamental theorem of the calculus}

 \begin{figure}[htbp] 
   \centering
   \includegraphics[width=4.in] {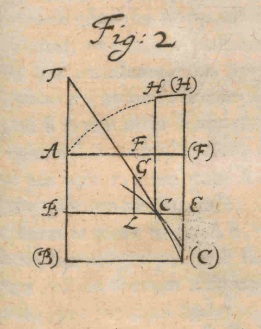} 
   \caption{Leibniz's  diagram  in the 1693 {\it Acta Eruditorum}. Courtesy of the Gottfried Wilhelm Leibniz 
   Bibliothek-Nieders\"{a}chsische Landesbibliothek- Hannover.}
   \label{fig:example}
\end{figure}

In this section I will contrast the geometrical proof   of the fundamental theorem of the calculus given by Leibniz (Leibniz, 1693), (Kowalewski, 1908), (Struik 1969,  282)
with the proof  of this theorem given  by Barrow  discussed in the previous section. 
Leibniz started his formulation of this theorem  with the statement
\beq
I shall now show {\it that the general problem of quadratures can be reduced to the finding of a line [curve] that has a given
law of tangency(declivitas)}, that is, for which the sides of the characteristic triangle have a given mutual relation.  Then
I shall show how this line [curve] can be described by a motion that I have invented. (Struik 1969, 282).
\enq
For his discussion\footnote{ In modern language, `` the general problem of quadratures"  is to obtain
the integral of a function $y(x)$,  described  by a curve with Cartesian coordinates $x,y$. By the fundamental
theorem of the calculus, this problem  can be solved
 by  finding another  function $z(x)$  described by a curve $x,z$,  with the property that its derivative $dz/dx=y$.
This property  corresponds to  Leibniz's  ``law of tangency" described by ``the characteristic triangle", 
which consist of an infinitesimal right angle triangle with  height to base ratio equal to $dz/dx$,
and hypotenuse aligned along the tangent of the curve $z(x)$.
}
   Leibniz introduced a diagram, Fig. 2,  that describes  two curves  
  $AH(H)$ and $C(C)$  with a  
common  abscissa  $AF$,  and   ordinates $HF$  and $FC$ respectively\footnote
{ There are two  errors, discussed in the Appendix,  in the reproduction of this  diagram
 in  Leibniz's mathematical papers edited by Gerhard (Leibniz 1693).
 In the corresponding diagram in Struik's {\it  A Source Book in Mathematics} (Struik, 1969), 
only one of these  errors  appeared.}.
Then Leibniz's  ``general problem of quadratures"  is  to obtain
the area bounded by  the curve  $AH(H)$  and the orthogonal lines  $AF$, $HF$,
by  finding the   curve $C(C)$  ``that has a given law of tangency,''  such that the  ordinate $FC$ of this curve  is proportional to 
this area.  Leibniz wrote that he will ``show how this line [$C(C)$] can be described by a motion that I have invented,"
and the  graphical device that he introduced for this purpose  will be described in the next section.\footnote
{In his  English translation of Leibniz's  1693 theorem, Struik  mentioned  that Leibniz described ``an instrument that can perform this construction", but  he
did not provide this description which is presented in Section 3.
}

The diagram associated with Leibniz's geometrical construction, shown in Fig. 2,  is essentially the same as  Barrow's diagram 
for  Proposition 11, Lecture 10,  shown in Fig. 1,  except that Leibniz's  diagram  is  rotated with respect to Barrow's 
diagram by $180^o$ around  the horizontal axis $AF$.\footnote
{ It might be thought that a geometrical proof 
of this   theorem would require the use of similar diagrams, but this is not necessarily the case. As a counter example,
compare Barrow's diagram, Fig. 1,  with  Newton's  diagram for his earliest geometrical proof of the fundamental theorem 
of the calculus (Guicciardini, 2009, 184).} 
Thus,  Leibniz's curves  $AH(H)$ and $C(C)$, Fig. 2,  correspond to Barrow's curves $ZEG'$, and  $AFI'$, respectively..
Likewise,  Leibniz's  tangent line $TGC$ to $C(C)$ corresponds to  Barrow's tangent line $TF$  to $AFI'$, see Fig. 3. 
Moreover, I will show that  Leibniz's  ``tangency law"  is the same as  the relation,  Eq. 3, applied by
Barrow in his proof that   $DT$ is the subtangent of $AFI'$.

 \begin{figure}[htbp] 
   \centering
   \includegraphics[width=4.in] {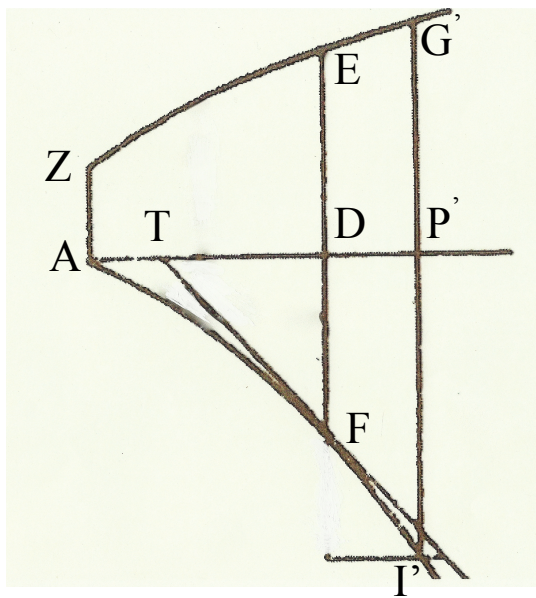} 
   \caption{Barrow's diagram, shown in Fig. 1, but   rotated by $180^0$ around the horizonal $APD$ axis with 
   some auxiliary lines deleted}
   \label{fig:example}
\end{figure}

Referring to his diagram, shown in Fig. 2, Leibniz continues
\beq
For this purpose  I assume for every curve C(C')  {\it a double characteristic triangle},  one $TBC$, that is assignable, and 
one, $GLC$, that is inassignable, and these two are similar.  The inassignable triangle consist of the parts $GL, LC$, with the elements
of the coordinates $CF, CB$ as sides, and $GC$ the element of arc, as the base of the hypotenuse. But the assignable triangle $TBC$
consists of the axis, the ordinate, and the tangent, and therefore contains the angle between the direction of the curve (or its tangent)
and the axis or base, that is, the inclination of the curve at a given point $C$. (Struik 1969, 283)
\enq

 The  tangent line to the curve $C(C)$ at $C$ is  
$TC$, and  $GL$ and $GC$ are  sides  of its ``characteristic triangle," $GLC$,  satisfying  a ``given mutual relation"
specified below. But  this  triangle which Leibniz called {\it inassignable} \footnote
{
 Leibniz chose the Latin  word {\it inassignabilis} for the characteristic triangle, because its
sides are differentials  which do not have an assignable  magnitude.
},
does not appear in the proof of Leibniz's theorem, while another characteristic triangle, $CE\bar C$, where $\bar C$
(not indicated in Leibniz's diagram, Fig. 2,  is
the intersection of the tangent line $TC$ with the extension of $(H)(F)$, turns out to be relevant to
Leibniz's proof.  Leibniz formulated  his  ``law of tangency'' 
in terms of the sides $TB$ and $BC$ of the similar  but  assignable triangle $TBC$.
The vertex   $T$ of this triangle is the intersection of the tangent $TC$  with a  line through the vertex $A$ 
of the curve $AH(H)$ perpendicular
to the abscissa $AF$, and $B$ is the intersection  with this line of a line from $C$ parallel to $AF$.
In Proposition 11, Lecture 10,  Barrow gave a proof of this law of tangency  which he  formulated, instead,  in terms of the
subtangent  $T'F$, where $T'$ (not shown in Fig. 2)  is the intersection of the tangent line $TC$
with the abscissa $AF$. 

Up to  this point,  $C(C)$ has been treated  as a {\it given} curve, but in the next sentence its {\it construction} is specified  by the
requirement that its slope  conforms to what Leibniz, in his introduction, called a certain ``law of tangency." 
\beq
Now let $F(H)$, the region of which the area has to be squared, be  enclosed between the curve $H(H)$, the parallel lines $FH$ and $(F)(H)$,
 and the axis $F(F)$, on that axis let  $A$ be a fixed point, and let a line $AB$, the conjugate axis, be drawn through $A$ perpendicular
 to $AF$. We assume that point $C$ lies on $HF$ (continued if neccesary); this gives a new curve $C(C)$ with the property  that,
  if from point  $C$ to the conjugate axis  $AB$ [an axis through $A$
perpendicular to $AF$] (continued if neccesary) both its ordinate $CB$ (equal to $AF$) and tangent $CT$ are
drawn, the part of the axis  between  them [$TB$] is to $BC$ as $HF$ to a constant $a$, or $a$ times $BT$ is equal
to the rectangle $AFH$ (circumscribed about the trilinear figure $AFHA)$. (Struik 1969, 283)
\enq
At this stage  the magnitude of $FC$,  the ordinate  of  $C(C)$,  was not specified, but further on  Leibniz announced   that $FC\cdot a$ is  the area of  the  
region $AHFA$,   where $a$ is an arbitrary constant with dimensions of length\footnote
{ The constant of proportionality  $a$ is the magnitude of
a fixed line which can be chosen arbitrarily, and  corresponds to  the chosen unit of length.}. 
The  relation  that the   curves $C(C)$ and $AH(H)$ must satisfy  is the requirement  
\be
\label{ratio1}
\frac{TB}{BC}=\frac{HF}{a},
\en
where $TB$ is the subtangent of $C(C)$ at $C$,  and $HF$ is the corresponding ordinate of the curve $AH(H)$ at $F$.
Presumably, this relation is Leibniz's `` law of tangency" that Leibniz announced in his introduction. 
Leibniz, however,  did not indicate  the origin of  this relation,  but setting $TB/BC=E(C)/EC$ where $E(C)C$
is the associated characteristic triangle, it  can be recognized 
as  the fundamental theorem of the calculus that Leibniz had obtained in differential form\footnote
{
In Leibniz's algebraic notation, setting $AF=BC=y$, $FC=x$, $EC=dy$ and $E(C)=dx$, by similarity of triangles $TBC$ and $CE(C)$, 
we have $TB/BC=E(C)/EC= dx/dy$. 
Then  setting $FH=z$, Eq. \ref {ratio1} becomes
$dx/dy=z/a$ which expresses the  fundamental theorem of the calculus in differential form.
}, 
and Barrow had proved\footnote
 {
 Formulated  in terms of the
subtangent  $T'F$, where $T'$ (a point  not labelled  in Leibniz's diagram)  is the intersection of the tangent line $TGC$
with the abscissa $AF$. By similarity of the triangles $TBC$ and $T'FC$ we have
$TB/BC=FC/T'F$,
which together with Leibniz law of tangency, Eq. 26,  gives the relation 
$T'F=a FC/HF$,
corresponding to Barrow's relation,  Eq. \ref{subt1},  with the parameter  $a$ replaced by $R$.
}in Proposition 11, Lecture 10.
But instead of giving  a proof of this relation, in the next sentence Leibniz just  asserted its validity:
\beq
This being established [the law of tangency, Eq. \ref{ratio1}], I claim that the rectangle on $a$ and $E(C)$  (we must discriminate between the ordinates $FC$
and $(F)(C)$ of the curve) is equal to the region $F(H)$. (Struik 1969, 254)
\enq

By ``the region $F(H)$'' Leibniz meant  the area  bounded by the
arc $H(H)$ of the curve $AFH$, the segment $F(F)$ of the abscissa, and the ordinates $HF$ and $(H)(F)$.
From the exact similarity of the  triangle $EC\bar C$ and the triangle $TBC$,  where $\bar C$ (not indicate in 
Leibniz's diagram, Fig. 2) is the intersection of the tangent line $TC$ with the
extension of the of the ordinate  $(H)(F)$, Fig. 2, it follows that 
\be
\label{char}
\frac{E\bar C}{EC}=\frac{TB}{BC},
\en
and  substituting in this relation Leibniz's form of the  law of tangency, Eq. \ref{ratio1},   yields
\be
\label{area1}
a\cdot E\bar C=EC\cdot HF.
\en
Since $CE=F(F)$, this expression  is the area of the rectangle inscribed in the region $F(H)$,  
which is smaller that the area of this region. Hence
\be
E\bar C < E (C),
\en
as is indicated in Leibniz's diagram,  Fig. 2,  but in the limit that $EC$ becomes vanishingly small,  the validity of Leibniz's 
assertion, quoted above, is established. 
Then the approximation that  $E(C)$ is the   differential change in the ordinate of $C(C)$ for a change $EC$ in the abscissa, 
and  $EC\cdot HF$, Eq. \ref{area1},   is the differential  area of the region $F(H)$, leads
to the differential form of the fundamental theorem of the calculus along the same lines
described by Barrow in Prop. 19, Lecture 11.   As Leibniz explained it,  

\beq
This follows immediately from our calculus.  Let $AF=y, FH=z, BT=t$, and $FC=x$; then
$t=zy/a$, according to our assumption [corresponding to Barrow's law of tangency, Eq.2,  in Leibniz's  coordinates]: on the other hand,
$t=ydx/dy$ because of the property of the tangents expressed in  our calculus. Hence $adx=zdy$ and therefore
$ax=\int z dy = AFHA$.(Struik 1969, 284)
\enq
Here $dx=E(C)$ and $dy=F(F)=CE$. Therefore,  $dx/dy=E(C)/EC$,  $t/y=BT/AF$,  and since
triangles $TBC$ and $(C)EC$ are similar, $TB/AF=E(C)/EC$,
which corresponds to $t=ydx/dy$.  Here Leibniz invoked this relation for the subtangent $t$  as a ``property of tangents expressed
in our calculus,"  but it  can be seen to  follow also from similarity relations between triangles in  his diagram.

In summary,  Leibniz based  his 1693  geometrical proof of the fundamental theorem of the calculus on  a ``law of tangency",  Eq. 11,
which  originally he had  developed in analytic form,  but  that also had been derived by Barrow  in his {\it Geometrical Lectures}, 
Prop. 11, Lecture 10. He then  proceed to demonstrate that a certain  curve 
$C(C)$, constructed according to this law,  gives the area bounded by a related   curve $AH(H)$, 
 along the same lines of Barrow's
proof of this  theorem in Prop. 19, Lecture 11,  using a diagram that turns out to be identical
 to Barrow's  diagram in Prop. 11, Lecture 10, after a  $180^0$ rotation around the horizontal axis (see Fig. 3)

\section {Leibniz graphical  device to perform integrations }
  \begin{figure}[htbp] 
   \centering
   \includegraphics[width=4in] {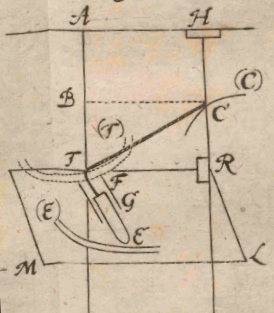} 
   \caption{ Leibniz's  diagram for a device to obtain graphically the area bounded by  a given curve  (1693 Acta Eruditorum). Courtesy of the Gottfried Wilhelm Leibniz 
   Bibliothek-Nieders\"{a}chsische Landesbibliothek- Hannover. }
   \label{fig:example}
\end{figure}

A novel feature of Leibniz's  1693 discussion of the fundamental theorem of the calculus is  his   description
of an ingenious  graphical  device to draw a  curve  $C(C)$ when  its slope  at each point is given. In his introduction, Leibniz mentioned that this curve ``can be described by a motion I have invented" \footnote
{ This description is missing in Struik's  English translation of Leibniz's 1693 geometrical proof of the fundamental theorem
of the calculus
}.
The motivation 
for this device appears to have been a problem that was brought to his attention  in 1676 by a physician, Claude Perrault,  namely,
to find the curve traced by a pocket watch pulled  along a  straight line  by its chain
 Leibniz solved the problem, but he did not divulge his solution until about 20 years later (Bos 1988, 9).

In this section we provide  an English translation of
Leibniz's description of his device, shown in Fig. 3,  obtained from the German  translation given in the
Oswald's Klassiker edition  of the exact sciences (Kowalewski, 1908,  32), 
and we also give a brief  discussion  of the principles  guiding its operation by an analytic description of  the curves involved.
\beq 
 In Fig. 4, the right angle $TAH$ is fixed and lies on a horizontal plane.  A vertical hollow cylinder  $TG$,  projecting out of  this plane, 
can move along the side $AH$.  On this cylinder  another massive cylinder  $FE$  slides upwards
and downwards  with a string $FTC$ attached at the  tip $F$, in such a way that  part $FT$ of the string lies  inside the hollow cylinder
and  part $TC$ lies on the above mentioned horizontal plane.  At the end $C$ of the string $TC$ there is a point (of a pen) that is lightly 
pressed against this plane to describes the curve $C(C)$. The movement  comes  from the hollow cylinder $TG$  that, while it
is guided from $A$ along $AT$ it  tightens  $C$.  The described point or pin $C$ pushes now $HR$ in the same horizontal plane
at right angle to $AH$ (the other side of the fixed right angle $TAH$) towards $A$  in a progressive manner.  This push does not
prevent that the point $C$ is moved only by the pull of the string, and therefore follows  its direction by this motion. There is also available a
board $RLM$,  that advances straight with the point $R$ perpendicularly  along  the staff $HR$,  after all continuously
driven by the hollow cylinder,  so that $ATHR$ is a rectangle.
Finally, there is a curve $E(E)$ described on this board (if you like in the form of a border)  in which the massive cylinder by means
of a cut  that one can image at the end $E$,  continuously intervenes; in this manner as $R$ moves towards $T$, the cylinder $FE$
moves upwards (along AT).  Since the length  $ET+TC$ is given (namely composed of the massive cylinder $EF$
and the entire string $FTC$), and  given the relation between $TC$ and $R$ or $BC$  (from  the given inclination [tangency]
law),  one obtains also the relation between $ET$and $TR$, the ordinate and abscissa of the curve $E(E)$  whose nature and
description on the board $LRM$ can be obtained through  ordinary geometry; one obtains also the description of the curve
$C(C)$  through this  available device. Now, however, it is  in the nature of our motion, that $TC$ always is tangent to the curve $C(C)$
so  the  curve $C(C)$ is described with the given  inclination law or the relation of the sides of the characteristic triangle
$TRC$ or $TBC$. Since this curve is the squaring figure corresponding  to the quadrature, as was shown a short while before,
one has obtained the desired quadrature or  measurement. (Kowalewski, 1908) 
\enq
Given the  slope of the curve $C(C)$, what is required to operate Leibniz's device is the relation between the ordinate
$ET$ and the abscissa $TR$ of the  curve $E(E)$. Leibniz states,  that this relation can be obtained  ``through ordinary geometry,"
leaving it as an exercise for the reader. Given $CR$ and $TR$, the  length $TC$
of the segment  of  the string lying on the plane $ATHR$ is
\be
TC=\sqrt{ TR^2+CR^2},
\en
and since the total length $TC+ET$ of the string has a  fixed value $U$,
\be
ET=U-\sqrt{ TR^2+CR^2}.
\en

The operation of  Leibniz's integration device can also be understood by expressing its variable components, string length $TC$
and inclination or slope $CR/TR$,  in analytic form. Introducing   Cartesian coordinates $x,y,z$ along the mutually
orthogonal axes $T_oR,T_oE,T_o A$ with origin at $T_o$, the initial location of the end point $T$ of the hollow cylinder,
the curve $C(C)$ is described by the coordinates $x=TR$,
$z=T_oB$, and the curve $E(E)$ by $x=TR$, $y=ET$. Then, by construction,
\be
\frac{dz}{dx}=\frac{CR}{TR},
\en
\be
TC=x\sqrt{1+(\frac{dz}{dx})^2},
\en
and
\be 
y=U- x\sqrt{1+(\frac{dz}{dx})^2}.
\en
For the special case that $ET$  is a constant,  
\be
\frac{dz}{dx}=-\frac{\sqrt{a^2-x^2}}{x}
\en
where $a=ET-U$, and the curve $C(C)$ is know as the $tractrix$.

\section {Concluding remarks} 
In his translation of the early mathematical manuscripts of Leibniz, J. M.  Child called attention to several  instances where 
Leibniz's  discussion and diagrams  are  similar to those  found in Barrow's {\it Geometrical Lectures} (Child 1920).  An example
is discussed in Appendix A. But 
 the relation  of Leibniz's  1693  {\it geometrical} proof  of the fundamental theorem of the calculus to the corresponding proof of Barrow
 has not been analyzed before.\footnote
{For example,  in his  translation of  Leibniz's   1693  article on the fundamental theorem of the calculus, D.J.  Struik wrote  that Leibniz's  ``expresses by  means of a figure the inverse relation of integration and differentiation " (Struik 1969, 282),   without commenting on the 
close similarity of this figure  to Barrow's  diagram in his   proof of this theorem (Struik 1969, 256).
}

In his  {\it Geometrical Lectures}, published in 1669,  Barrow gave two  elegant geometrical proofs of the fundamental theorem of the calculus
 discussed here in detail.
In Proposition 11, Lecture 10, he established this theorem by  rigorous geometrical bounds,  while 
in Proposition 19, Lecture 11, he gave an alternate proof  based on the application of  infinitesimals along the lines
practiced by mathematicians in the 17-th century.  The eminent historian of science  D.T.  Whiteside, however,  
belittled Barrow's  proofs   calling them   ``the work of a competent
university don" (Whiteside 1961, 289),  but  this evaluation is contradicted  by  the evidence presented here.  
 Likewise, another historian of science, M. S.  Mahoney,  concluded that Barrow ``was not particularly original"  (Mahoney 1990, 180,240),
without, however,  providing
any evidence that  Barrow's  proofs  had been given earlier.  Actually, Barrow was
among the foremost mathematicians of his time, whose misfortune was to  be eclipsed by his former prot\'{e}g\'{e},  Isaac Newton.
 In contrast to Whiteside's  and Mahoney's  remarks, a distinguished mathematician,
Otto Toeplitz, concluded that ``in a very large measure Barrow  is indeed the real discoverer  [of the fundamental theorem of the calculus] - insofar as an individual can ever be given credit within a course of development such as we have tried to trace here" (Toeplitz  1963, 98).
 After   reviewing the historical
record,  M. Feingold   also disagreed with Whiteside's  and Mahoney's  views (Feingold, 1990).  
Recently,  Guicciardini   discussed 
some of  Barrow's geometrical  proofs  showing that these proofs also  had had a greater impact on Newton's own contributions to the development
of the calculus  than had been realized in the past (Guicciardini, 2009).  But Toeplitz also commented that ``what Newton absorbed from
the beginning remained foreign to Barrow throught his life:  the turn from the geometrical to the computational function concept -
the turn from the confines of the Greek art of proof to the easy flexibility of the indivisibles. On one page of his work Barrow
alluded briefly to these matters, but quickly, as though in horror, he dropped  them again.(Toeplitz 1963, 130)
 
Twenty three  years after Barrow published his work, Leibniz  presented a  very similar 
geometrical proof of the fundamental theorem of the calculus
 theorem  based on a  ``law of tangency" that  Leibniz gave as ``being established".   But  this law  corresponds 
 to a theorem that Barrow had proved in Proposition 11, Lecture 10.
 Moreover, as has been  shown here, Leibniz's  diagram, Fig. 2,  is essentially  the same, apart from
 orientation\footnote
 { An anonymous referee of an earlier version of this manuscript  suggested  that the similarity  between Leibniz's and Barrow's diagram follows from a common tradition, 
 originating with H. van  Heuraet (Heuraet 1637),  to represent geometrically   the area bounded by one curve by a second curve.  But
 such a  diagram  can be drawn  in many different ways as can be seen, for example, in one of Newton's diagram,  which is 
 quite different from Barrow's,  illustrating  his first  geometrical proof of the fundamental theorem of the calculus.
 (Guicciardini 2008,  184). I  thank N. Guicciardini for first calling van Heuraet's diagram to my attention.
 }
as Barrow's diagram,  Fig.1, given in  this proposition,  see also Fig. 3,  
and Leibniz's  arguments, which were  based on differential quantities, are  the same as  those given by 
 Barrow  in  Prop. 19, Lecture 11. 
 Since Leibniz had obtained a copy of Barrow's {\it Geometrical Lecture}  in 1673,  it is  implausible that in the intervening twenty years he had
 never encountered  Barrow's two propositions, and the fact  that his  diagram and geometrical proof of the fundamental theorem of the calculus
 are virtually  the same as Barrow's is unlikely to be  a  coincidence.  In fact,  Leibniz's marginal annotations in his copy
 of Barrow's book indicate that at least by 1676 he had studied one of Barrow's propositions contained in one of his
 last lectures,  Prop. 1, Lecture 11 (Child 1920,  16) (Leibniz 2008, 301),  see Appendix A.
 Child wrote that ``as far as the actual invention of the calculus as he understood the
term is concerned, Leibniz received no help from Newton or Barrow; but
for the ideas that underlay it, he obtained from Barrow a great deal more
than he acknowledged, and a very great deal less than he would have
like to have got, or in fact would have got if only he would have been
more fond of the geometry he disliked. For, although the Leibnizian
calculus was at the time of this essay far superior to that of Barrow
on the question of useful application, it was far inferior in the
matter of completness." (Child 1920, 136). Also  Feingold commented
that  ``I find it difficult to accept that historians
can argue categorically that books that
a person owned  for years went unread simply because
he or she failed failed to find
dated  notes from these books - especially when
the figure in question is Leibniz, who
was truly a voracious reader. And how can one determine
with certainty what a genius
like Leibniz  was capable of comprehending from various
books and letters he encounter
or discussion he participated in, however confused their context
appears to us today?
Such reasoning, it seems to me, substitutes preconceived notion
for constructive  historical knowledge" (Feingold  1990, 331).
Finally,  it   should  not be forgotten that 
  Leibniz's, when composing  his own version
  of planetary motion,   {\it Tentamen de motuum coelestium causis} (Nauenberg 2010, 281),  also denied having read Newton's {\it Principia}, but his denial has been shown to be false (Bertoloni-Meli,  1993).
 Hence, Leibniz's persistent claim,  particularly  during his priority controversy with Newton,   of not having  any indebtedness 
 to Barrow in his development of the calculus,  must be taken {\it cum grano salis}.
 
 To his credit,  in his 1693 article in the {\it  Acta Eruditorum},  Leibniz also called attention to the usefulness of the fundamental theorem of the calculus  for the  evaluation of  integrals, and for this purpose he designed  a device to evaluate  integrals graphically,
 see Section 4. Moreover, his work stimulated  the applications  
 of the calculus by  his celebrated contemporaries,
the Bernoulli brothers,  Jacob Herman, and Pierre Varignon to the solution of problems in mechanics (Nauenberg, 2009).
In his own development  of the fundamental theorem of the calculus, Newton also 
realized the great usefulness of this  theorem for integration, and for this purpose he created  extensive
   tables of integrals. But  he kept  these results  to himself, and he did not publish them
   until 1704 when he appended  them
  to his $\it Opticks$  his  {\it Two Treatises on the Species and Magnitudes of Curvilinear Figures}  (Whiteside 1981, 131)

Finally,  it should be emphasized  that until the 19th century,  when
 the {\it analytic}  calculus was established on proper mathematical foundations\footnote{
 Leibniz never was able to give a proper definition of his differentials which were derided by  Bishop  George Berkeley  as  the 
``ghost of departed quantities" ( Boyer 1989).  For a discussion of the importance of proper foundations for the analytic form
of the calculus,  and its  development in the nineteenth century, see  (Grabiner 1981).
 }, Barrow's Prop. 10 Book 2 was already
 a  rigorous proof, based on sound {\it geometrical}  principles,  of its  fundamental theorem\footnote{  
 At about the same time, a  similar proof  was given by James Gregory (Baron, 1969, 233).}. Today, however, 
 partly due to the dismissive remarks about Barrow by  historians of science like  Whiteside and Mahoney, (Whiteside 1961;Mahoney, 1990), it is Leibniz and Newton who get most of  the credit for the development of
 the calculus  while Barrow has been more or less forgotten.  But, to  quote Rosenberger,
 \beq
  Like all great advances in the sciences,
 the analysis of the infinitesimals 
 did not suddenly arise,  like Pallas Athena out of the head of Zeus,  from the genius of a single author, but instead  it 
 was carefully prepared and  slowly grown, and finally after  laborious trials  by the strength of genius,  its general significance
 and long range meaning was brought to light" (Folkerts 2001, 299).
\enq

\begin{figure}[htbp] 
   \centering
   \includegraphics[width=5.in] {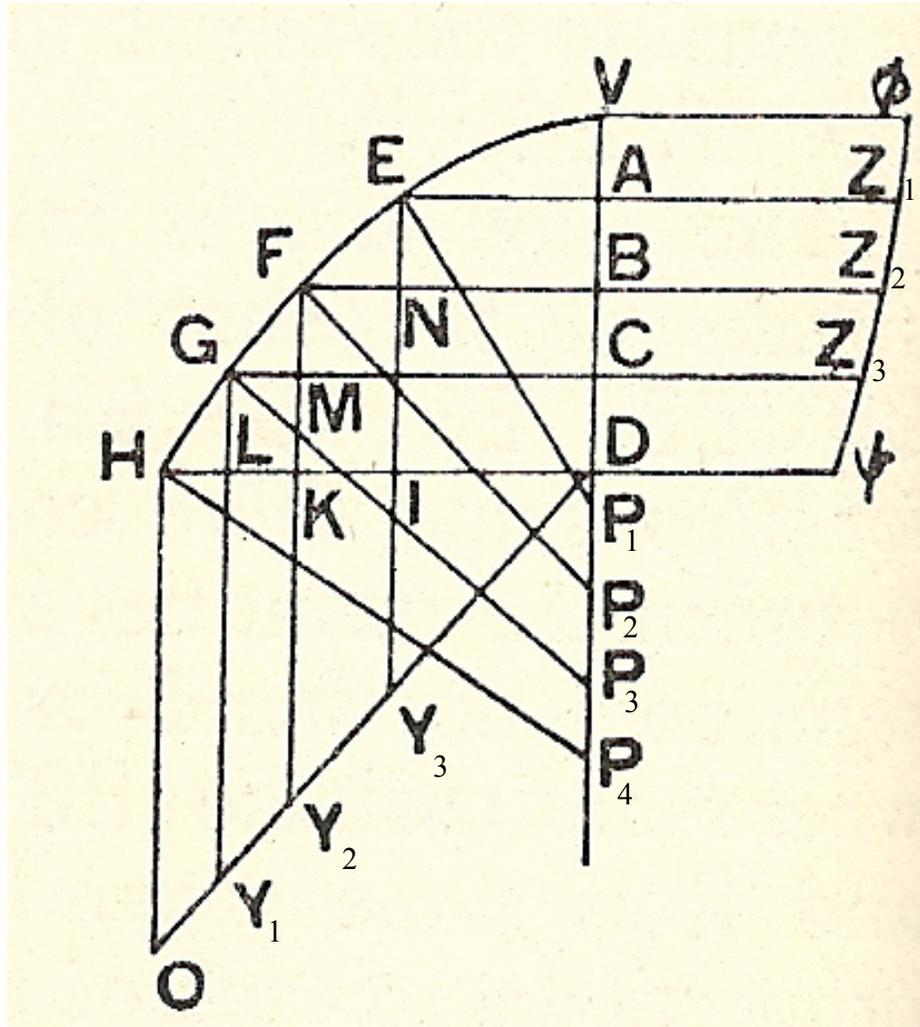} 
   \caption{ Barrow's diagram for the integration of a  curve of subnormals }
   \label{fig:example}
\end{figure}

\section*{Appendix A. Area of a curve of subnormals} 
In Proposition 1, Lecture 11,  Barrow presented an ingenious geometrical construction to obtain  the area bounded by
 a  curve  $\phi Z \psi$ with
ordinates equal to the  subnormals of  a {\it given} curve $VEH$,  with common abscissa $VD$,  shown in Fig. 5.
For clarity, we  have added numerical subscript  to  the  symbols $P, Z$ and $Y$ that appear repeatedly in Barrow's  original diagram.
Barrow gave a  proof  that this area is
equal to the area  of a  right angle isosceles triangle $HDO$ with sides equal to $HD$  shown at the bottom 
of his diagram, where $HD$ is the largest ordinate of $VEH$.  The horizontal and vertical lines  in Barrow's diagram,  Fig.5, 
illustrate his  approximation by rectangles of the required areas,
 while  the diagonal lines are the normals to $VEA$ at the chosen values of the  abscissa  of this curve.
  For example,  at $A$  the subnormal $AP_1$ is obtained by finding the intersection at $P_1$  of the normal to  $VEH$ at $E$ with the axis $VD$,  where $EA$ is the ordinate $A$. The corresponding ordinate of 
$\phi Z \psi$  is constructed  by setting $AZ_1=AP_1$, where $AZ_1$ is taken along the extension of $EA$.
When $E$  approaches  $V$, the arc $EV$ can be approximated
by a straight line, and the ``characteristic"   triangle $VAE$, assumed to be infinitesimal, becomes  similar to the triangle $EAP$.
This construction is then repeated at the  equally space points  $B,C$ and $D$ along the abscissa $VD$, e.g.
at the next point $B$ on the abscissa, the normal is $BP_2$,  where this second value of $P$ is obtained by the intersection of the normal to $VED$ at $F$ with $VD$.  In this case the
characteristic triangle is $ENF$  which is similar to $FBP$,  and $BZ_2=BP_2$.
Then, in Barrow's words,
\beq
the space $VD\psi \phi$ differs in the least degree only from the sum of the rectangles $DC\cdot D\psi+CB \cdot CZ_3+BA \cdot  BZ_2 +AV \cdot AZ_1$
\enq

Setting the intervals $AV,BA,CB,$ and $DC$ that appear to be equal as  $\Delta x$,  Barrow's  approximation
to the area $VD\psi \phi$ is given by the area of  the  sum  of rectangles, $\Delta x \cdot ( AP_1+BP_2+CP_3+DP_4)$.   Then, in the limit that          $\Delta x$ becomes vanishingly small, and the number of
rectangles increases indefinitely, this sum becomes equal to  the   area $VD \psi \phi$
The similarity
of the characteristic triangles with the triangles associated with the subnormals implies that  $EA/VA=AP_1/EA, FN/EN=BP_2/FB,
GM/FM=CP_3/CG$ and $HL/GL=DP_4/HD$.  Hence, the sum $\Delta x \cdot ( AP_1+BP_2+CP_3+DP_4)= EA\cdot EA+FN\cdot FB+
GM\cdot GC+HL\cdot HD$. 
According to Barrow's construction, Fig. 5, the intervals  $EA=ID,  FN=KI, GM=LK$  and $HL$, are unequal,
and  $LY_1=GC, KY_2=FB$ and $IY_3=EA$. Hence, the above sum  is  equal to $HL\cdot HO+
LK\cdot LY_1+KI \cdot KY_2 +ID\cdot IY_3$
which corresponds to a sum of rectangles, giving  an upper bound to the area of the right angle isosceles triangle $HDO$.
In the limit that the  interval $\Delta x$ become vanishingly small, and the number of rectangles that bound both
areas increases indefinitely
this sum gives the area of the triangle $HDO$, $(1/2) HD^2$,  leading to Barrow's conclusion that
\be
area\; VD \psi \phi     = (1/2)HD^2
\en

By always expressing   position of points on his diagrams by letters,  sometimes repeating the same letter for different points,
Barrow lacked  a suitable notation to describe his sums,  particularly in  the  limit $n\rightarrow \infty$.  Moreover, for this reason
the  relations that he had obtained geometrically  between
two finite sums was not evident algebraically \footnote {By setting $x_i, y_i$ for  the abscissa and ordinate of a point,  $\Delta x_i=x_{i+1}-x_i$,
$\Delta y_i=y_{i+1}-y_i$, and
$n_i=y_i \cdot (\Delta y_i)/\Delta x_i) $ for the subnormal at this point,  where $i=1,n$,  the relation between Barrow's two sums 
becomes obvious without his geometrical analysis:
\be 
 \sum^{i=n}_{i=1} \Delta x \cdot n_i = \sum ^{1=n}_{i=1} \Delta x_i\cdot y_i \cdot (\Delta y_i)/\Delta x_i=\sum^{i=n}_{i=1} \Delta y_i  \cdot y_i ,
 \en
where
\be
 \sum^{i=n}_{i=1} \Delta x \cdot n_i \approx  area \,VD \psi \phi,
 \en
 and
 \be
 \sum^{i=n}_{i=1} \Delta y_i  \cdot y_i \approx (1/2) HD^2,
 \en
 In this
proposition Barrow only hinted at the limit $n\rightarrow \infty$  with the intriguing remark,
\beq
A lengthier indirect argument may be used  but what advantage is there?
\enq
But in an appendix
to lecture 12,  he  discussed more carefully  the upper and lower  bound of  the area of a curve, referring to an
{\it indefinite} number of rectangles. Later on, Isaac Newton improved  Barrow's discussion,  and included it as   Lemma 2, Book 1,  in the
{\it Principia}  (Guicciardini 2009, 178, 221)}

Although it lacked the rigorous mathematical justification of Barrow's analysis,
the power of Leibniz's  useful notation is that by the substitution of his relation
 $dx\cdot  n(x)= dy \cdot y$ between the differentials $dx$ and $dy$,  it reduces Barrow's lengthy  geometrical construction
and derivation to a one line  analytic relation between two integrals (Leibniz 1686)
\be
\int^x_0 dx' \; n(x')=\int^y_0 dy'\; y'= \frac{1}{2}y^2.
\en

We have shown  that Leibniz's relation
is based on the  characteristic triangle discussed by Barrow,  that gave rise to the  relation
$n_i/y_i \approx \Delta y_i/\Delta x$,  for $i=1,n$ .
Leibniz claimed  that he first  learned  about  the  characteristic triangle from  Pascal, but  evidently he must
have recognized it also when he examined  Barrow's diagram  (Child 1920, 16).

Like Barrow,   Leibniz also  labelled  points on his geometrical diagram with letters,  but in the case that the same letter
appeared repeated, he added a number of parenthesis corresponding to the number of times
this letter was repeated, e.g. $C,(C), ((C))$ etc.  But  later, he also distinguished repeated letters by adding a 
numerical subscript  in front of these letter,  e.g. $_1C, \;_2C, \;_3C$  (Child 1920,137; Leibniz 2008, 573; Bertoloni Meli 1993, 109,135)
\begin{figure}[htbp] 
   \centering
   \includegraphics[width=3.in] {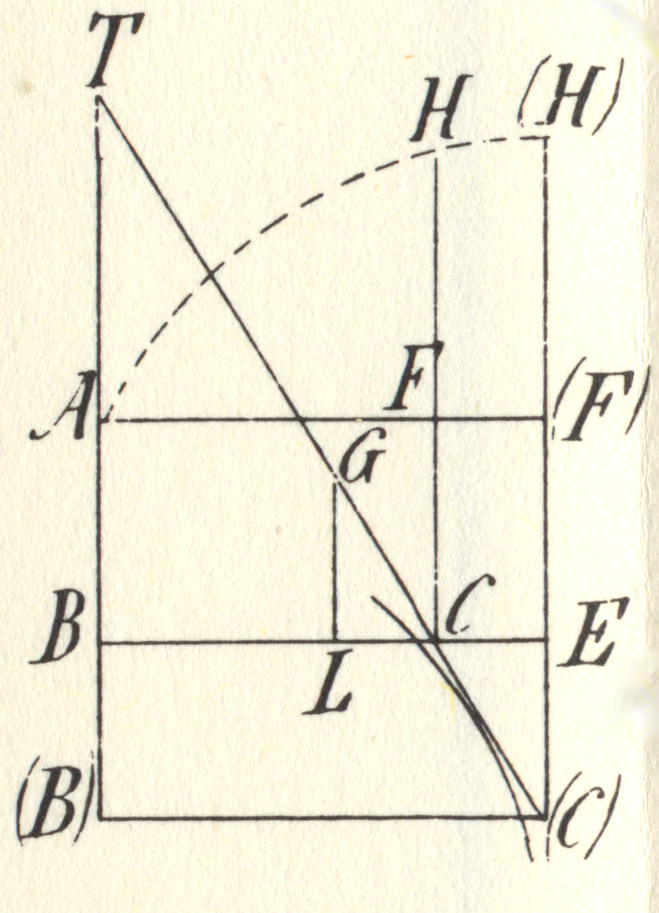} 
   \caption{ Reproduction of Leibniz's  1693 diagram  in  Gerhardt's edition of Leibniz's mathematical papers (Leibniz, 1693)  }
   \label{fig:example}
\end{figure}

  \begin{figure}[htbp] 
   \centering
   \includegraphics[width=3.in] {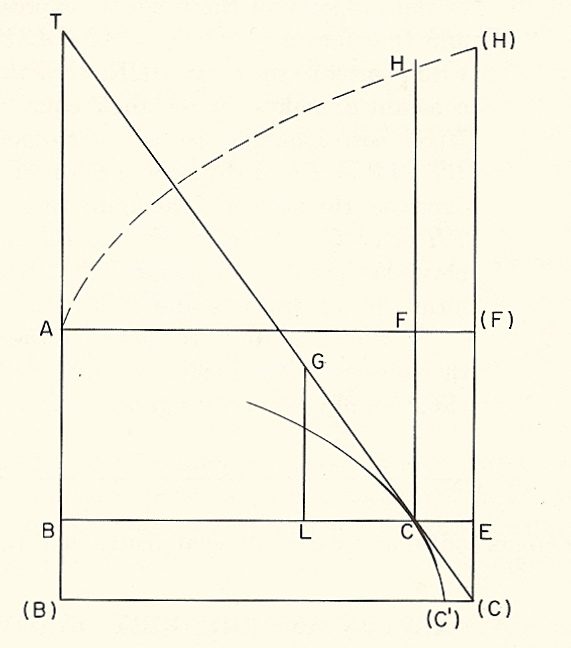} 
   \caption{ Reproduction of Leibniz's  1693 diagram  in (Struik 1696, 283) }
   \label{fig:example}
\end{figure}

 \begin{figure}[htbp] 
   \centering
   \includegraphics[width=3.in] {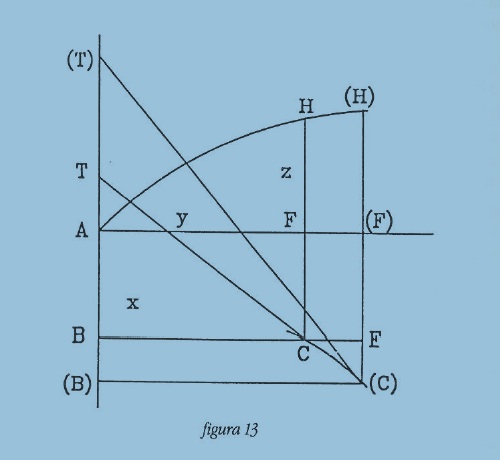} 
   \caption{ Reproduction of Leibniz's  1693 diagram  in (Giacardi 1995, 322) }
   \label{fig:example}
\end{figure}

\section* {Appendix B. Errors in the reproduction of Leibniz's diagram} 
It should be pointed out  that Leibniz 1693  diagram, Fig. 2,  for his geometrical proof of the fundamental theorem 
of the calculus,  has  been  repeatedly reproduced incorrectly.  

In Gerhardt's edition of Leibniz's mathematical
papers, this reproduction, shown in  Fig. 5,  contains two errors: 1)  the curve labelled $C(C) $ 
touches  tangentially the line $TC$  below the intersection
$C$  of  this line with $BE$, and 2) the extension of $TC$  ends at  its intersection  
with the extension of $(F)(H)$,  labelled incorrectly $(C)$. But $(C)$ is  the end point of the curve $C(C)$, and the line  $E(C)$ is
proportional to the the area $FH(H)(F)$ which is greater than the distance between $E$ and the intersection
of the extension of $TC$ that was  labelled $\bar C$ in Section 2.

In Struik's reproduction (Struik 1969, 283), shown in Fig. 6, Leibniz's curve labelled $C(C')$
is drawn  correctly,  but  the intersection $(C)$  is again shown incorrectly 
 as the extension of the tangent  $TC$ intersecting the extension,
 of $(F)(H)$.  Thus,  when Struik  translates Leibniz's text 
 \beq
 This being established, I claim that the rectangle of $a$ and $E(C)$ . . . is equal to the region $F(H)$.
 \enq  
 it appears as if  Leibniz had made a mistake here, but this is due  to reference to  Struik's incorrect diagram.
 
 More  recently, L. Giacardi  also reproduce Leibniz's diagram incorrectly, see Fig 8,  drawing  the extension of the
 tangent line $TC$ of $C(C)$  at $C$ to intersect this curve at $(C)$, and by also introducing a line $T(C)$, supposedly tangent
 to $C(C)$ at $(C)$,  which does not even appear
 in Leibniz's diagram or in his text (Giacardi 1995, 322)
  
 Such  errors make  Leibniz's text difficult to comprehend. 

\section{Acknowlegdments}
I would like to thank Niccol\`{o} Guicciardini for many helpful comments and discussions;  Eberhard Knobloch  and
two anonymous  referees  for helpful  criticisms of an earlier version of this manuscript, and  for 
providing some relevant  references. I also thank Mordechai Feingold for  some helpful suggestions.

\end{document}